\documentclass{article}
\usepackage[utf8]{inputenc}
\usepackage[english]{babel}
\usepackage{csquotes}
\usepackage{biblatex}
\usepackage{hyperref}

\usepackage{dirtytalk}

\usepackage{ dsfont }
\usepackage{amssymb}
\usepackage{amsmath}        % rozšíření pro sazbu matematiky
\usepackage{amsfonts}       % matematické fonty
\usepackage{amsthm}         % sazba vět, definic apod.
\usepackage{bbding}         % balíček s nejrůznějšími symboly
\usepackage{bm}             % tučné symboly (příkaz \bm)
\usepackage{graphicx}       % vkládání obrázků
\usepackage{fancyvrb}       % vylepšené prostředí pro strojové písmo
\usepackage{indentfirst}    % zavede odsazení 1. odstavce kapitoly
\usepackage{icomma}         % inteligetní čárka v matematickém módu
\usepackage{dcolumn}        % lepší zarovnání sloupců v tabulkách
\usepackage{booktabs}       % lepší vodorovné linky v tabulkách
\usepackage{paralist}       % lepší enumerate a itemize
\usepackage{xcolor}         % barevná sazba
\usepackage{ dsfont }

\input xy                         %pro diagramy
\xyoption{all}
\usepackage{amsmath,amscd}
\usepackage{textcomp}
\theoremstyle{plain}
\newtheorem{veta}{Věta}
\newtheorem{Thm}[veta]{Theorem}
\newtheorem{Prop}[veta]{Proposition}
\newtheorem{Ex}[veta]{Example}
\newtheorem{Cor}[veta]{Corollary}
\newtheorem{Lemma}[veta]{Lemma}

\theoremstyle{plain}
\newtheorem{Def}[veta]{Definition}
\newtheorem{Remark}{Remark}

\theoremstyle{remark}

\theoremstyle{plain}

\newcommand{\C}{\mathbb{C}}

\newcommand{\Q}{\mathbb{Q}}
\newcommand{\Z}{\mathbb{Z}}

\newcommand{\m}{\mathfrak{m}}

\newenvironment{dukaz}{
  \par\smallskip\noindent
  \textit{Proof}.
}{

\rightline{$\qedsymbol$}
}

\usepackage{blindtext}
\title{Formal matrix representations of rings with a Nakayama permutation}
\author{Dominik Krasula
}

\begin{document}
 \maketitle 
\begin{abstract}
Semiperfect rings with Nakayama permutations, which include pseudo-Frobenius and Frobenius rings in particular, are studied by applying the general theory of formal matrix rings to their Peirce decompositions. We utilise this approach to derive a novel combinatorial criterion that determines whether such a ring has a prescribed Nakayama permutation and, subsequently,  whether it has an essential right socle. 

Such a ring can be represented as a block matrix ring, where the blocks on the diagonal correspond to the cycles in the Nakayama permutation. We characterise all possible supports for such blocks. Faithful pairing bimodules with essential simple socles appear on the shifted diagonal of their formal matrix representations, thereby generalising the known description of pseudo-Frobenius rings. 

We introduce a novel method for glueing rings with a Nakayama permutation, a construction based on our combinatorial criterion.  This construction is used to create an indecomposable Frobenius ring with two simple modules, whose endomorphism rings are not isomorphic, resolving a question left open in our previous work.
\end{abstract}
\smallskip


\smallskip

\noindent \textbf{Author:} RNDr. Dominik Krasula

Charles University, Faculty of Mathematics and Physics, 

Department of Algebra

Sokolovská 49/83, 186 75 Praha 8, Czech Republic

krasula@karlin.mff.cuni.cz

ORCID: 0000-0002-1021-7364

\smallskip 

\noindent This work is a part of project SVV-2025-260837.

\noindent This research was supported by the grant GA ČR 23-05148S from the Czech
Science Foundation.

\smallskip 

%\noindent \textbf{Declarations of interest}: I have nothing to declare. 

\section{Introduction}

T. Nakayama, in his seminal article [15],  defined \textit{quasi-Frobenius rings}, or QF rings for short, as artinian rings with a Nakayama permutation. QF rings are characterised as self-injective artinian rings. It is the class of rings in which injective and projective modules coincide.

A ring $R$ is \textit{semiperfect} if and only if the unity $1_R$ can be written as a finite sum of pairwise orthogonal local (hence primitive) idempotents. A permutation $\pi$ on this set is called a \textit{Nakayama permutation} if for any local idempotents $e$ and $f:=\pi(e)$, the pair $(eR, Rf)$ is an \textit{i-pair}, that is, $soc(eR)\cong top(fR)$ and $soc(Rf)\cong top (Re)$. 

This decomposition of $1_R$ induces a representation of a semiperfect ring as a \textit{formal matrix ring} (also called a \textit{ring of generalised matrices}) with local rings on the diagonal (Secs. \ref{SecSemiperfect},  \ref{SecMatrix}). This idea originates from the work of B. J. Müller [13, Sec. 2].

The core contribution of this article is a combinatorial criterion (Thm.  \ref{MainThm}) to determine whether a formal matrix ring has a prescribed Nakayama permutation. This is proven using existing theory on formal matrix rings, as presented in [10]. 

 \medskip

Our research aims to gain a deeper understanding of the information that can be inferred about the ring from its Nakayama permutation. Pioneering work in this area was done by K. R. Fuller in [4], whose work led to the development of the concept of i-pair. See [1] for related results and literature on this topic. Building on Fuller's results, A. T. Hannula proved that corners of basic QF rings corresponding to a cycle in a Nakayama permutation are again basic QF and described how to glue them together by means of Morita dualities.

W. Xue generalised Fuller's results to semiperfect rings with essential socles [24]. If such a ring has a Nakayama permutation, it is self-injective if and only if it is linearly compact [24, Thm. 7], and this is equivalent to a \textit{pairing bimodule} $eRf$ being linearly compact for each pair of orthogonal primitive idempotents $e, f\in R$.

The class of semiperfect rings with Nakayama permutation and essential socles strictly contains \textit{Pseudo-Frobenius} rings, or PF rings for short (also called \textit{cogenerator rings}) [12, Thm. (19.25)],  and \textit{dual rings} (also called \textit{annihilator rings}), see [6, Thm 4.2, Lemma 3.2].

Using the properties of Morita duality, Müller provided a full description of PF rings whose Nakayama permutation is a cycle in [14, Thm. 6]. He generalised and extended Hannula's results for PF rings in [14]  using B. Roux's generalisation of Fuller's results, [19]. In a PF ring $R$, $(eR, Rf)$ is an i-pair if and only if $eRf$ defines Morita duality between the local corner rings $eRe$ and $fRf$.  We refer the reader to Xue's monograph [23] for terminology and basic results on Morita duality. A ring with duality is semiperfect [23, Cor. 3.14, Thm. 4.3] and formal matrix rings with a duality are described in [23, Thm. 4.12].

\medskip 

We continue this study by introducing a combinatorial criterion rather than the use of Morita dualities. The main combinatorial criterion is presented in Theorem \ref{ThmCombEss}, which checks whether a ring with a Nakayama permutation has an essential right socle.

Our combinatorial approach generalises Müller's results (Sec.  \ref{SecInjCycle}) and provides a straightforward method of glueing rings with a Nakayama permutation together (Prop. \ref{PropGlue}), cf. [7, Prop. 3.3, Thm. 3.6] and [14, Section 3]. We answer the question posed in Remarks after [14, Thm. 6] in Corollary \ref{CorGlue}.

If $R$ is a ring with a Nakayama permutation and essential socles, the corners induced by cycles in the permutation are themselves rings with a Nakayama permutation (Thm. \ref{ThmCycle}). Note that in such rings, the left and right socle coincide (Remark \ref{RemarkTT'}).

Motivated by this, rings whose Nakayama permutations are cycles are studied further.  It is shown that any local ring with a duality can be a local corner of a ring with a Nakayama permutation in Example \ref{ExCycle}. Pairing bimodules on the shifted diagonal are shown to be faithful with simple coincident essential socles (Remark \ref{RemarkTT'} and Corollary \ref{CorFaithful}). A complete characterisation of possible supports of these rings is given in Section \ref{SecInjCycle}. If the socles coincide, simple modules linked by a Nakayama permutation have isomorphic endomorphism rings (Prop. \ref{PropK}).

 Our method of glueing rings with a Nakayama permutation (Sec. \ref{SecGlue}) implies that the structure of corner rings corresponding to different cycles of the Nakayama permutation can be quite different, as long as the endomorphism rings of their simple modules are compatible (Thm. \ref{ThmGlue}). It follows that indecomposable rings with arbitrary Nakayama permutations exist in various classes of rings, such as \textit{semiprimary} but not QF, \textit{right perfect} but not semiprimary, etc.  
 
 An example of a QF ring such that the endomorphism rings of their simple modules are not all isomorphic is constructed (Ex. \ref{Exnon-isomfields}), thus answering a question left open in the author's previous work [9, Prop. 11].

\medskip

The article is organised as follows. Section \ref{secFormal} sets the formal framework for the article.  Section \ref{SecQF} provides an easy combinatorial criterion  (Thm. \ref{MainThm}) characterising whether a formal matrix ring has a prescribed Nakayama permutation and some of its immediate applications. Section \ref{SecEss} refines the combinatorial criterion for rings with essential socles (Thm. \ref{ThmCombEss}).  

Section \ref{SecGlue} introduces a general constructive method for glueing  Nakayama permutations to a new indecomposable ring that respects their respective permutations. The main result is Theorem \ref{ThmGlue}, but the construction is more general (Cor. \ref{CorGlue}).

\section{Preliminaries and notation} \label{secFormal}

This section defines formal matrix rings and sets the notation used throughout the text. The primary reference is [10, Chapter 2.3], but we restrict our definition of formal matrix rings, as we always assume that the rings on the diagonal are local. This is in contrast to the general definition given in [10], which we refer to as \textit{block matrix rings}. A ring has such a representation if and only if it is a semiperfect ring.

Subsection \ref{SecSemiperfect} explains the representations of semiperfect rings as formal matrix rings, and Subsection \ref{BasicSec} shows that we can restrict our attention to basic rings.   The formal matrix rings and modules over them are defined in Subsections \ref{SecMatrix} and \ref{SecModules}.

When referring to a side-specific condition (such as artinian, injective, etc.) of a ring or a bimodule without specifying a side, we always mean \textit{as both left and right modules}. We use [11, Chapters 7, 8] as a reference for the properties of idempotents and semiperfect rings.

\subsection{Formal matrix representations of semiperfect rings}\label{SecSemiperfect}

Recall that a unital ring $R$ is semiperfect if and only if $1_R$ can be written as a sum of finitely many pairwise orthogonal local idempotents. In such rings, an idempotent is local if and only if it is primitive [13, Thm. 1].

Given a complete set of pairwise orthogonal local idempotents $e_1,\dots, e_m$ in $R$, we have the following series of isomorphisms:
\[R\cong End_R(R_R)=Hom(\bigoplus_{1 \leq i \leq m } e_iR, \bigoplus_{1 \leq i \leq m } e_iR)\cong \bigoplus_{1\leq i, j \leq m} Hom_R(e_jR, e_iR)\cong  \bigoplus_{1\leq i, j \leq m} e_iRe_j \]
in the category of abelian groups, sometimes referred to as a \textit{Peirce decomposition of $R$}.

Using this decomposition, we can represent $R$ as a set of square matrices of order $m$ such that for $i, j\leq m$, an element in the $(i,j)$-th entry is an element of an $R_i\text-R_j$-bimodule $e_iRe_j$, where $R_i:=e_iRe_i$ is a local ring. The standard matrix addition and multiplication correspond to the addition and multiplication in the formal matrix representation of $R$. Multiplication  $e_ire_j\cdot e_jr'e_k=e_ire_jr'e_k$ induces an $R_i\text-R_k$-bimodule homomorphism $e_iRe_j\otimes_{R_j}e_jRe_k\to e_iRe_k$. Given a right $R$-module $M$, we have an isomorphism $M\cong Me_1\oplus \dots \oplus Me_m$ in the category of abelian groups, and $Me_i$ is a right $R_i$-module via the restriction of the right $R$-action on $M$.

This representation is not unique, as the decomposition of $1_R$ is determined up to reordering and conjugation by a unit.

\subsection{Basic rings}\label{BasicSec}

This section recalls the notion of a basic ring for a semiperfect ring. Any categorical property of a ring $R$ can be tested on its basic ring. One can construct a formal matrix representation of any ring Morita equivalent to $R$ from that of its basic ring. The key ring properties studied in this article are categorical: semiperfect with a Nakayama permutation $\pi$, right self-injective, PF, QF, possessing a right (self-)duality [23, Cor. 4.6], right perfect and right artinian (Proposition \ref{PropStructure}). However, combinatorial properties may differ. All rings Morita equivalent to a Frobenius ring $R$ are Frobenius if and only if $R$ is \textit{weakly symmetric}, i.e., its Nakayama permutation is the identity. 

\smallskip

Let $R$ be a semiperfect ring, and $T_1,\dots, T_n$ be a list of all simple right $R$-modules up to isomorphism. Then we can choose an ordering of a complete set of pairwise orthogonal idempotents $e_1,\dots, e_m$ in such a way that the projective cover of $T_i$ is isomorphic to $e_iR$ for any $i\leq n$. Then $e_1R,\dots, e_nR$ is a complete set of finitely generated indecomposable projective right $R$-modules up to isomorphism.  Two right $R$-modules $e_iR$ and $e_jR$ are isomorphic if and only if $Re_i\cong Re_j$ as left $R$-modules.

The idempotents $e_1,\dots, e_n$ are then called a \textit{basic set of idempotents}. To this set, we assign \textit{multiplicities}, a set of positive integers $\mu_1,\dots, \mu_n$ such that 
 \begin{gather*}
     R_R\cong (e_1R)^{\mu_1}\oplus \dots \oplus  (e_nR)^{\mu_n} \label{ModuleDecomp} \\ 
          _R R\cong (Re_1)^{\mu_1}\oplus \dots \oplus  (Re_n)^{\mu_n},
 \end{gather*}
in the category of right and left modules, respectively.

The ring $End_R(e_1R\oplus \dots \oplus e_nR)$ is isomorphic to $eRe$, where $e=e_1+\dots + e_n$, and is called a \textit{basic ring for $R$}. Any ring Morita equivalent to $R$ is isomorphic to $End((e_{1}R)^{m_1}\oplus \dots \oplus (e_{n}R)^{m_n})$ for some positive integers $m_1,\dots, m_n$ [12, Prop (18.33)]. Up to isomorphism, Morita equivalent rings have the same set of local corner rings.

\subsection{Formal matrix rings} \label{SecMatrix}

A ring $R$ is \textit{local} if and only if its \textit{Jacobson radical} $J(R)$ equals the set of all non-units. A simple local ring is then a \textit{skew field}. For a local ring $R$, we assign a \textit{residue (skew) field} $\m$, together with a canonical projection $R\twoheadrightarrow R/J(R)=:\m$. This projection induces an $\m$-module structure on $R$-modules.  
\begin{Def}\label{Defmain}
Let $n$ be a positive integer and $R_1,\dots, R_n$ be a set of local rings with residue skew fields $\m_1,\dots, \m_n$, respectively. For $i, j, k\leq n$, let $B_{i,j}$ be a \emph{coordinate $R_i\text-R_j$-bimodule} and \[\varphi_{i,j,k}\colon B_{i,j}\otimes_{R_j} B_{j,k} \to B_{i,k}\] an $R_i\text-R_k$-bimodule \emph{product homomorphism.} We set $B_{i,i}:=R_i$ and $\varphi_{i,i,j}$ ($\varphi_{i,j,j}$) is induced by the left $R_i$ (right $R_j$) action on bimodule $B_{i,j}$.

Furthermore, assume that the following \emph{associativity condition}
\[(\alpha) ~~~~\varphi_{i,k,l}(\varphi_{i,j,k}(a\otimes_{R_j} b) \otimes_{R_k}c) =   \varphi_{i,j,l}(a\otimes_{R_j} \varphi_{j,k,l}(b\otimes_{R_k}c))\]
holds for any $i, j,  k, l\leq n$ and any $a\in B_{i,j}$, $b\in B_{j,k}$ and $c\in B_{k,l}$.

Then we define a \emph{formal matrix ring $R$ of order $n$} as a set of square matrices of order $n$ such that for $i, j\leq n$ an element in the $(i,j)$-th entry is an element of  $B_{i,j}$.  Addition and multiplication are analogous to the standard matrix addition and multiplication. The product homomorphisms facilitate the multiplication of elements from different bimodules.

For $i, j\leq n$, we denote by $E_{i,j}$ a function mapping an element $a\in B_{i,j}$ to an element in $R$ such that all positions are zeros, except that the $(i,j)$-th position is $a$.

We assign to $R$ a complete set of pairwise orthogonal local \emph{canonical idempotents} $e_1,\dots, e_n$, where $e_i:=E_{i,i}1_{R_i}$.
\end{Def}
To simplify the notation, we write $ab$ in lieu of $\varphi_{i,j,k}(a\otimes_{R_j}b)$.Observe that $B_{i,j}$ is isomorphic to $e_iRe_j$ as an $R_i\text-R_j$-bimodule. Thus, a formal matrix ring is its formal matrix representation. Note that using the tensor product is merely a compact way to ensure that the distributivity axioms hold, and we do not work with $B_{i,j}B_{j,k}$ as we do with a tensor product, as we are always interested only in pure tensor elements. 

\begin{Remark} \label{AssRemark}
    In condition $(\alpha)$, it is enough to assume that $i\neq  j \neq k \neq l$; otherwise, it follows from the associativity of the ring actions and the axioms of the tensor product.
\end{Remark}
We can construct a formal matrix ring from any finite set of local rings and corresponding bimodules by setting all \textit{proper} product homomorphisms $\varphi_{i,j,k}$ with $i\neq j \neq k$ to zero. We call such rings \textit{trivial formal matrix rings}.  These correspond to semiperfect rings $R$ such that for any triple of local idempotents $e, f, g\in R$, such that $ef=0=fg$, the equality $eRfRg=0$ holds.  All Frobenius trivial formal matrix rings are characterised in Corollary \ref{TrivCor}.
\begin{Ex}
    The ring $M_n(S)$ of square matrices of order $n$ over a local ring $S$ is an example of a formal matrix ring.  The product homomorphisms are induced by the multiplication in $S$ and the condition ($\alpha$) is then automatic. Note that, except for the case $n=1$, these are not basic rings, so the combinatorial results of this article will not apply. 
\end{Ex}
Being (right) artinian, semiprimary, and (right)
perfect are examples of ring-theoretic properties that are independent of the choice of product homomorphisms, see Section \ref{SecStructure}. Another example is rings with a duality, see [23, Thm. 4.12].

\smallskip

From a formal matrix representation of a ring, we can deduce its decomposition in the category of rings. An idempotent $e\in R$ is \textit{central}, if it commutes with all elements in $R$. This is equivalent to $eR(1-e)=0=(1-e)Re$, which means that the corresponding entries in the formal matrix are zero. Consequently, $R$ decomposes as $eRe\oplus (1-e)R(1-e)$ in the category of rings.

\subsection{Modules over formal rings}\label{SecModules}

 Following [5, Thm. 1.5], we can represent any right $R$-module $M$ as a row  module $(M_1,\dots, M_n)$, where $M_i$ is a right $R_i$-module, together with a right $R_j$-homomorphism \[f_{i,j}\colon M_i\otimes_{R_i} B_{i,j} \to M_j\]   for each $ i \neq  j \leq n$ called \textit{proper homomorphism of module multiplication}, such that the ring action is associative, i.e.,
\[(\beta)~~~~~ f_{j,k}(f_{i,j}(m\otimes_{R_i} b)\otimes_{R_j} c)=f_{i,k}(m\otimes_{R_i} bc) \]
 holds for any  $i\neq j\neq k \leq n$ and any $m\in M_i$, $b\in B_{i,j}$ and $ c\in B_{j,k}$.  An $R_i$-homomorphism of module multiplication $f_{i,i}\colon M_i\otimes_{R_i} R_i\to M_i$ is induced by the right $R_i$-action on $M_i$.  The left modules are represented as column modules in an analogous manner. We write $mb$ in lieu of $f_{i,j}(m\otimes_{R_i} b)$.

By $\epsilon_{i}$, we denote the standard basis vector where the $i$-th entry is $1_{R_i}$ and the remaining entries are zero. We can view it as a functor from $Mod\text-R_i$ to $Mod\text-R$ mapping an $R_i$-module $M_i$ to an $R$-module $\epsilon_i M_i$, i.e., a row module whose entries are all zero except for the $i$-th entry, which is $M_i$.
 
Let $N=(N_1,\dots, N_n)$ be a right $R$-module. Then $N$ is an $R$-submodule of $M$ if and only if for any $i, j\leq n$ it holds that $N_iB_{i,j}$ is an $R_j$-submodule of $M_j$. The case $i=j$ amounts to $N_i$ being an $R_i$-submodule of $M_i$.

An $R$-homomorphism $\chi\colon X\to M$ is a collection of $R_i$-homomorphisms $\chi_i\colon X_i\to M_i$ compatible with proper homomorphisms of module multiplication in the following sense
\[(\gamma)~~~~~ \chi_j(xb)=\chi_i(x)b,\]
for any $i\neq j\leq n$ and any $x\in X_i$ and $b\in B_{i,j}$.

Recall that in a category of modules, a homomorphism is a monomorphism (isomorphism) if and only if it is injective (bijective). It follows that a homomorphism $\chi$ is mono (iso) if and only if each component $\chi_i$ is. 

\section{Combinatorial criterion}\label{SecQF}

This section proves Theorem \ref{MainThm}, giving a combinatorial criterion characterising whether a basic formal matrix ring has a prescribed Nakayama permutation. Key observations about socles are in Subsection \ref{SecSocle}. Artinian and perfect formal matrix rings are described in Subsection \ref{SecStructure}. 

Throughout the text, the symmetric group on the set $\{1,2,\dots, n\}$ is denoted by $S_n$.

\begin{Def}\label{DefQF}
    Let $R$ be a semiperfect ring and $e_1,\dots, e_n$ be a basic set of idempotents with multiplicities $\mu_1,\dots, \mu_n$ and, let $\pi\in S_n$.

    Then $R$ is a \emph{ring with a Nakayama permutation $\pi$} if $soc(e_iR)\cong top(e_{\pi(i)}R)$ and $soc(Re_{\pi(i)})\cong top(Re_{i})$ for each $i\leq n$.

    A ring $R$ is \emph{QF} if and only if it is an artinian ring with a Nakayama permutation. Furthermore, if $\mu_{i}=\mu_{\pi(i)}$ for each $i\leq n$, $R$ is called a \emph{Frobenius} ring.
\end{Def}
Nakayama permutation is unique up to conjugation. A right artinian ring with a Nakayama permutation is artinian [24, Cor. 11] and hence QF. For examples of rings with Nakayama permutations that are not QF, see Example \ref{ExCycle}. The following remark compiles some well-known properties of QF rings used repeatedly throughout the text.

    \begin{Remark}\label{RemarkQF}    
Because all multiplicities are 1 in a basic ring, it follows that a basic QF ring is Frobenius.
    
          If $R$ is a QF ring, its left and right socles coincide [12, Cor. (15.7)]. $soc(R)$ is, as an $R\text-R$-bimodule, isomorphic to $top(R)=R/J(R)$ if and only if $R$ is Frobenius. In particular, a local artinian ring is Frobenius if and only if its right (left) socle is simple as a right (left) ideal.
    \end{Remark}
Idempotents $e_i$ and $e_{\pi(i)}$ can be viewed as \textit{linked} by the Nakayama permutation. Due to repeated use, we propose the following coinage. Do not confuse this with the concept of linked idempotents from the block decomposition theory [11, Sec. 22].
\begin{Def}
    Let $R$ be a basic formal matrix ring of order $n$ with a Nakayama permutation $\pi$ and let $i,j \leq n$.

    We then call $B_{i,\pi(i)}$ a \emph{pairing coordinate bimodule} and $\varphi_{i,j,\pi(i)}$ a \emph{pairing product homomorphism}.
\end{Def}

\subsection{Socles}\label{SecSocle}

This subsection presents several straightforward observations about simple modules over a basic formal matrix ring, thereby providing a description of the socles of modules generated by primitive idempotents used in the proof of Theorem \ref{MainThm}.

\smallskip
 
Let $R$ be a semiperfect ring. Then $R$ is semilocal, i.e., its top is totally decomposable as a ring.  By the Wedderburn–Artin Theorem, if  $R$ is basic with a basic set of idempotents $e_1, \dots, e_n$, then there exist skew fields $\m_1,\dots, \m_n$ such that
\[top(R)\cong top(e_1Re_1)\oplus \dots \oplus top(e_nRe_n)\cong \m_1\oplus \dots \m_n,\]
 in the category of rings. Representing $R$ as a formal matrix ring, we have
\[top(R)\cong ({\m_1}, \dots, \m_n)\cong \bigoplus_{1\leq i\leq n}\epsilon_i   \m_i  \]
in the category of right $R$-modules. Thus, we proved the following.
\begin{Lemma}\label{lemmaSimple}
    Let $R$ be a basic formal matrix ring of order $n$ and let $T$ be a simple right $R$-module. 

    Then there exists $i\leq n$ such that $T$ is isomorphic to a module $\epsilon_i \m_i$ with all proper homomorphisms of module multiplication being zero.  
\end{Lemma}
 The left version is analogous. A simple left module is isomorphic to $\epsilon^T_i  \m_i$ for some $i\leq n$.
 \begin{Remark}\label{BasicRemark}
     The assumption that $R$ is basic is necessary. Otherwise, there is $i\leq n$ such that $\mu_i>1$ and $M_{\mu_i}(\m_i)$  is isomorphic to a direct summand of $top(R)$ in the category of rings. Then, the corresponding simple right $R$-module is of the form $\epsilon_{i_1}\m_i + \dots + \epsilon_{i_{\mu(i)}}\m_i$ for some set of distinct indices $i_1,\dots, i_{\mu(i)}\leq n$.

     In particular, if $D$ is a skew field and $n>1$, then all simple right $M_n(D)$-modules are isomorphic to $(D, D, \dots, D)$.
 \end{Remark}
We can now formulate a crucial step in the proof of Theorem \ref{MainThm}.
\begin{Lemma}\label{Lemasocle}
Let $R$ be a basic formal matrix ring, and let $M=(M_1,\dots, M_n)$ be a right $R$-module.

For each $i\leq n$, let $N_i$ be the $R_i$-submodule of $M_i$ consisting of elements $m\in M_i$ such that $mB_{i,j}=0$ for all $j\neq i$. 

Then $soc(M)=(soc(N_1), \dots, soc(N_n))$.
\end{Lemma}
\begin{dukaz}
Module $M$ decomposes as $M=\bigoplus_{1\leq i \leq n} \epsilon_i M_i$. By the previous lemma, it follows that a simple $R$-submodule of $M$ is of the form $\epsilon_i T$ for some $i\leq n$, where $T$ is a simple $R_i$-submodule of $M_i$. Thus, $T\subseteq N_i$, because $\epsilon_i T$ is closed under the right $R$-action. The conclusion then follows from the fact that simple submodules of $M$ generate $soc(M)$.
\end{dukaz}
We conclude this section with the following observation, which characterises basic formal matrix rings, cf. [14, Lemma 3], and is central to the proof of Theorem \ref{ThmGlue}.
\begin{Lemma}\label{LemmaJacob}
    Let $R$ be a basic formal matrix ring of order $n$ and $i\neq j \leq n$.

    Then $B_{i,j}B_{j,i}=Im(\varphi(i,j,i))$ is contained in $J(R_i)$.
\end{Lemma}
\begin{dukaz}
From the description of $top(R)$, we see that $J(R)$ can be described as matrices whose elements on the diagonal are from the Jacobson radicals of the corresponding local rings, cf. [10, Thm. 2.4.3]. In particular, matrices with zeros on the diagonal are in $J(R)$.

Let $m\in B_{i,j}$ and $n\in B_{j,i}$. We can view them as elements in $R$ by representing them as matrices $E_{i,j}m$ and $E_{j,i}n$ with only one nonzero entry. Then $E_{i,j}m\in J(R)$ and $E_{j,i}n\in J(R)$ and hence  $E_{i,i}mn$ is contained in $J(R)$, showing $mn\in J(R_i)$.
\end{dukaz}

\subsection{Structural properties of formal matrix ring} \label{SecStructure}

In this short detour, we characterise when a formal matrix ring is artinian or perfect.  The artinian case is proved in [10, Thm. 2.5.1], whereas the perfect case is proved only for trivial formal matrix rings [10, Cor. 3.7.13] despite being seemingly well known, see [13, Thm. 2], cf. [18, Props. 2.7.35, 2.7.14(iii)].
\begin{Prop}\label{PropStructure}
    Let $R$ be a formal matrix ring of order $n$.
    
    Then $R$ is right artinian (noetherian) if and only if all rings $R_i$ are right artinian (noetherian) and for each $i, j\leq n$, the bimodule $B_{i,j}$ is finitely generated as a right $R_j$-module.

    The ring $R$  is semiprimary (right perfect) if and only if all rings $R_i$ are semiprimary (right perfect).
\end{Prop}
Note that a right perfect self-injective ring is QF [16, Thm. 6.39]. The question of whether a semiprimary right self-injective ring is QF is an open problem known as \textit{Faith's Conjecture}. For a discussion and literature review, see [3]. A right perfect ring satisfies the descending chain condition for finitely generated left ideals [18, Thm. 2.7.36]. Therefore, a right perfect left noetherian ring is left artinian.

\subsection{Main theorem}\label{SecCriterion}

This section formulates a combinatorial criterion to test whether a formal matrix ring $R$ has a prescribed Nakayama permutation and fixes the notation used throughout the text. Condition (2) describes properties of pairing coordinate bimodules, whereas property (1) describes the remaining coordinate bimodules.

\begin{Thm}[Combinatorial criterion] \label{MainThm} Let $R$ be a basic formal matrix ring of order $n$ and let $\pi\in S_n$.

Then $R$ is a ring with a Nakayama permutation $\pi$ if and only if for each $i \leq n$: 

(1) If $j\neq \pi(i)$, then either 

(i) for any nonzero element $m\in B_{i,j}$ there exists an index $ k\neq j \leq n$ such that $mB_{j,k}\neq 0$ or

(ii) the right $R_j$-submodule of $B_{i,j}$ consisting of all elements $m$ such that $mB_{j,k}=0$ for all $k\neq j$ has zero socle.

(1')  If $j\neq \pi(i)$, then either 

(i')  for any nonzero element $m\in B_{i,j}$ there exists an index $h\neq i \leq n$ such that $B_{h,i}m\neq 0$ or

(ii') the left $R_i$-submodule of $B_{i,j}$ consisting of all elements $m$ such that $B_{h,i}m=0$ for all $h\neq i$ has zero socle.

(2) The right socle $T_{\pi(i)}$ of the $R_i\text-R_{\pi(i)}$-subbimodule $N_{\pi(i)}$ of $B_{i,{\pi(i)}}$ consisting of elements $m$ such that $mB_{\pi(i), j}=0$ for all $j\neq \pi(i)$ is simple as a right $R_{\pi(i)}$-module.

(2') The left socle $T'_i$ of the $R_i\text-R_{\pi(i)}$-subbimodule $N'_{i}$ of $B_{i, \pi(i)}$ consisting of all elements $m$ such that $B_{h,i}m=0$ for all $h\neq i$ is simple as a left $R_i$-module.
\end{Thm}
Evidently, condition (ii) implies (i). This formulation of the theorem was chosen to stress the importance of condition (i). We are primarily interested in rings with an essential right socle. For such rings, condition (i) and (i') always hold with $k:=\pi(i)$ and $h:=\pi^{-1}(j)$; see Theorem \ref{ThmCombEss}.
\begin{dukaz}  
    We prove Conditions (1) and (2) in the forward implication. The rest is analogous.
    
    Assume that $R$ has a Nakayama permutation $\pi$ and let $i\leq n$. Then $e_iR$ has a simple socle isomoprhic to $top(e_{\pi(i)})R$. Because $R$ is basic, it follows that there exists a simple right $R_{\pi(i)}$-submodule $T\leq B_{i,\pi(i)}$ such that $\epsilon_{\pi(i)}T$ is the socle of $e_iR=(B_{i,1}, \dots, B_{i,n})$. Part (2) thus follows directly from Lemma \ref{Lemasocle}. 
    
    Suppose that (1) does not hold, i.e., for some $j\neq \pi(i)$  subbimodule $N_j$ of $B_{i,j}$ consisting of all elements $m$ such that $mB_{j,k}=0$ for all $k\neq j$ has a nonzero socle. Then $\epsilon_j soc(N_j)$ is a nonzero semisimple right $R$-submodule of $e_iR$. This is a contradiction with $soc(e_iR)\cong top(e_{\pi(i)}R)$.
\end{dukaz}
%In the next section, we present a version of the combinatorial criterion for the case when socles are essential. Müller has a deeper description of PF rings whose Nakayama permutation is a cycle in [BJM], some of which we recall in the next section.

We end this section with an immediate, though
somewhat surprising, result that will be recalled in Section \ref{SecGlue}. The case for PF rings was observed by Müller [14, Thm. 6, Remark].
\begin{Prop}\label{PropFixedPoint}
    Let $R$ be a basic formal matrix of order $n$ with a Nakayama permutation $\pi\in S_n$ and let $i\leq n$ be an index such that $\pi(i)=i$.

    If $R_i$ is simple, then it is a direct summand of $R$ in the category of rings.
\end{Prop}
\begin{dukaz}
Because $i=\pi(i)$ and $R_i$ is simple, it follows that $N_i=N_{\pi(i)}$ equals the whole $R_i$ by the condition (2). But all modules over a simple ring are faithful, so it implies that $B_{i,j}=0$ for any index $i\neq j\leq n$. Similarly, by (2') we have that $B_{j,i}=0$. Thus, the canonical idempotent $e_i$ is central and $R$ decomposes as $e_iRe_i\times (1-e_i)R(1-e_i)$ in the category of rings.
\end{dukaz}
A result by I. Reiten [17, pp. 72-73] implies that an indecomposable QF artin algebra has a simple injective module if and only if it is a skew field. The following generalises this.
\begin{Cor}
    Let $R$ be a ring with a Nakayama permutation, and let $S$ be a simple right injective module. 

    Then there exists a central idempotent $e\in R$ such that $S\cong eRe$.
\end{Cor}
\begin{dukaz}
We begin by proving that there exists an idempotent $e$ such that $S \cong eR$.  
    Without loss of generality, we assume that $R$ is basic. By definition \ref{DefQF}, it follows that $R$ is a right Kasch ring. So, there is an injection $S\hookrightarrow R_R$. Let $S'$ be the image of this injection.  Because $S'\leq R$ is injective, it is a direct summand of $R$. Furthermore, it is a right faithful ideal of $R$. Thus, it is not nilpotent. This means that $S'$ is a simple idempotent ideal of $R$, so it is generated by a local idempotent.  Let us denote this idempotent $e$.

    We find a basic set of idempotents in $R$ such that $e_1=e$ and consider a formal matrix representation of $R$ with respect to these idempotents, and let $\pi$ be the corresponding Nakayama permutation. Then $\pi(1)=1$ and $e_1Re_1$ is simple. Thus, $e_1=e$ is central. 
\end{dukaz}
We conclude this section with the following class of examples. Recall that if a ring $S$ has a self-duality defined by an $S\text-S$-bimodule $E$, then $E$ is an injective cogenerator [23, Thm. 2.4(3)]. In particular, if $S$ is local, $E$ can be chosen to be the injective envelope of its unique simple module.
\begin{Ex}\label{ExCycle}
Let $S$ be a local ring that is commutative or has a self-duality, and let $E$ be an injective envelope of its unique simple module.

Then the trivial formal matrix ring $R:=\begin{pmatrix}
 S& E & 0 & 0 & \dots \\
 0 & S & E & 0 &\dots \\
 \vdots & & \ddots & \\
 E &\dots &&0&S    
\end{pmatrix}$ of order $n$ is a ring  with a Nakayama permutation $\pi=(1~2~\dots~n)$.    
\end{Ex}
The two conditions on $S$ coincide when it is linearly compact [23, Thm. 6.8]. If $S$ is artinian and commutative, then $E$ is a finitely generated module [12, Thm. (3.64)], thus by Proposition \ref{PropStructure}, $R$ is artinian and hence Frobenius. 

Recalling Proposition \ref{PropStructure}, if $S$ is semiprimary but not artinian, we obtain an example of a semiprimary ring with a Nakayama permutation that is not QF. Similarly, we can construct a perfect ring with a Nakayama permutation that is not semiprimary.

 \section{Essential socle}\label{SecEss}

 In this section, we focus on rings with a Nakayama permutation whose right socle is essential. To this end, we formulate a combinatorial criterion, characterising them by means of pairing coordinate bimodules and pairing product homomorphisms. This criterion demonstrates that many properties of PF rings, as proven in [14], can be generalised to this setting. The properties established here will serve as the assumptions for the main theorem in the next section.
 
\subsection{Combinatorial criterion} \label{SecEssComb}

Recall the notation from Theorem \ref{MainThm}.
\begin{Thm}\label{ThmCombEss} Let $R$ be a basic formal matrix ring of order $n$ with a Nakayama permutation $\pi$.

Then the right socle of $R$ is essential as a right ideal if and only if for each $i\leq n$

(a) The simple right $R_{\pi(i)}$-module $T_{\pi(i)}$ is an essential right $R_{\pi(i)}$-submodule of $B_{i,\pi(i)}$.   

(b) For all $j\neq \pi(i)$ and nonzero $m\in B_{i,j}$, it holds that $mB_{j,\pi(i)}\neq 0$.    

The left socle of $R$ is essential as a left ideal if and only if for each $i\leq n$

(a') The simple left $R_i$-module $T'_{i}$ is an essential left $R_i$-submodule of $B_{i,\pi(i)}$.   

(b') For all $j\neq i$ and all nonzero $m\in B_{j,\pi(i)}$, it holds that $B_{i,j}m\neq 0$.    
\end{Thm}
It follows from (a) that $T_{\pi(i)}$ is a unique simple right $R_{\pi(i)}$-submodule of $B_{i,\pi(i)}$ for each $i\leq n$.
It follows from (b) that for a ring with a right essential socle, Condition (1) from a combinatorial criterion is satisfied by means of Condition (i).
\begin{dukaz}
We only prove the \say{right} part of the theorem. Firstly, the right socle of $R$ is essential as a right ideal if and only if the socle of each $e_iR$ is essential. Indeed, any nonzero ideal $I\leq e_iR$ must contain a simple submodule of $e_iR$ if the socle of $R$ is essential. For the converse, observe that if $I\leq R$ is a nonzero right ideal, then there exists $i\leq n$ such that $e_iI$ is a nonzero right submodule of $e_iR$.

Thus, it is enough to show that  $T_{\pi(i)}$ is an essential submodule of $e_iR$ for each $i\leq n$. We claim that this is equivalent to conditions (a) and (b). Recall that a submodule is essential if and only if it has a nonzero intersection with all nonzero cyclic submodules. Furthermore, it is enough to consider elements with only one nonzero coordinate, because each $\epsilon_jm_j$ is contained in $mR$, as shown by multiplying $m=(m_1,\dots, m_n)$ by $E_{j,j}$.

Condition (a) is a necessary condition. If there is $b\in B_{i,\pi(i)}$ such that $bR_{\pi(i)}\cap T_{\pi(i)}=0$, then $bR$ also has zero intersection with $T_{\pi(i)}$. Hence, $bR$ has zero intersection with the entire right socle, as each element $bR$ lies in a coordinate bimodule $B_{i,j}$ for some $j$, but $B_{i, \pi(i)}$ is the only such coordinate bimodule that contains elements of the socle. 

We now prove that, assuming (a), condition (b) is equivalent to $T_{\pi(i)}$ being an essential submodule. Take $j\neq \pi(i)$ and a nonzero $m\in B_{i,j}$. Assuming (a), the module $mR$ has zero intersection with $T_{\pi(i)}$ if and only if it has zero intersection with $B_{i,\pi(i)}$. Equivalently, $mR\cap B_{i,\pi(i)}=mB_{j,\pi(i)}=0$.
\end{dukaz}
\begin{Remark}\label{RemarkTT'}
    If the ring has a Nakayama permutation and an essential right socle, its socles coincide (see the proof of [NY, Thm. 3.21]).

By [2, Thm. 2.5], if such a ring is right noetherian, it is right artinian and hence a QF ring. 
\end{Remark}
The following was known in the PF case [14, Thm. 4] as the pairing coordinate bimodules $B_{i,\pi(i)}$ induce a duality between the corresponding local corners. Note that \textit{right} essentiality implies \textit{left} faithfulness.
\begin{Cor}\label{CorFaithful}
     Let $R$ be a basic formal matrix ring of order $n$ with a Nakayama permutation $\pi$ and let $i\leq n$. Further, assume that its right socle is essential as a right ideal. 

     Then $B_{i,\pi(i)}$ is a faithful left $R_i$-module for any $i\leq n$ such that $i\neq \pi(i)$.
\end{Cor}
\begin{dukaz}
If $i=\pi(i)$, then $B_{i,\pi(i)}=R_i$ which is clearly faithful as a module over itself. For $i\neq \pi(i)$, choose $j=i$ in the condition (b) of the combinatorial criterion.
\end{dukaz}
In this setting, pairing coordinate bimodules are faithful and indecomposable on both sides, and they have essential simple socles that coincide. However, they need not be injective [6, Ex. 6.1, 6.2] , so they do not need to be cogenerators. 

In the case of rings with an essential right socle, Theorem \ref{ThmCombEss} provides a full description of trivial formal matrix rings with a Nakayama permutation.
\begin{Cor}\label{TrivCor}
    Let $R$ be an indecomposable, basic, trivial formal ring of order $n>1$ with an essential right socle and a Nakayama permutation $\pi$. 

    Then $\pi$ is an $n$-cycle, and the pairing coordinate bimodules are the only nonzero proper coordinate bimodules. 
\end{Cor}
\begin{dukaz}
    Let $i\neq j\leq n$ such that $j\neq \pi(i)$. If  $m\in B_{i,j}$ is nonzero, then by condition (b), we have $mB_{j,\pi(i)}\neq 0$. But $R$ is a trivial matrix ring, thus $B_{i,j}$ cannot contain any nonzero element. 
    
    Now let $e\in R$ be the sum of all idempotents in a cycle of a $\pi$. By the previous observation, it is a central idempotent. If $R$ is indecomposable, this means that $e$ is conjugate to $1_R$. So the Nakayama permutation is a cycle.
\end{dukaz}
As a straightforward corollary of Theorem \ref{ThmCombEss}, we see that an indecomposable ring with a Nakayama permutation cannot be represented as a triangular matrix ring. 
\begin{Prop}\label{PropTriangular}
    Let $R$ be a basic formal matrix ring with an essential right socle of order $n$ with a Nakayama permutation $\pi$. Set $e:=e_1+\dots + e_k$ for some $k\leq n$.

    Then $(1-e)Re=0$ if and only if $eR(1-e)=0$.
\end{Prop}
\begin{dukaz}
    We prove the forward implication by contraposition. The proof of the other implication is symmetric. Assume that there exists $i\leq k < j$ such that $B_{i,j}$ is nonzero. From the essentiality of the socle and Condition (b), we have $B_{i,j}B_{j,\pi(i)}\neq 0$ and, in particular, $B_{j,\pi(i)}\neq 0$.
    
   If $\pi(i)\leq k$, the proof is complete. Otherwise, $\pi(i)>k$. Now consider a cycle in $\pi$ containing $i$. It contains numbers less than or equal to $k$ and greater than $k$. Hence, there is some $l$ such that $\pi^l(i)>k\geq \pi^{l+1}$ and thus the corresponding pairing coordinate bimodule $B_{\pi^l(i), \pi(\pi^l(i))}$ is a nonzero module in the bottom left block $(1-e)Re$.
\end{dukaz}

\subsection{Rings whose Nakayama permutation is a cycle}\label{SecInjCycle}
We apply Theorem \ref{ThmCombEss} to show that rings whose Nakayama permutations are cycles are the building blocks of rings with Nakayama permutations, and we further investigate their structure.  

\begin{Thm}\label{ThmCycle}
    Let $R$ be a basic formal matrix ring of order $n$ with a Nakayama permutation $\pi$ and essential socles. Furthermore, let $I\subseteq \{1,2,\dots, n\}$ be a nonempty set of indices such that if $i\in I$ then $\pi(i)\in I$ and set $e:=\sum_{i\in I}e_i$. 

    Then $eRe$ is a ring with a Nakayama permutation $\pi_{\restriction_I}$ and essential socles.

    Furthermore, if $R$ is PF/QF, then so is $eRe$.
\end{Thm}
For QF rings, this was observed in [7, Prop. 2.1] based on  [4, Lemma 2.2]. The case for PF rings was proven in [14, Lemma 5].
\begin{dukaz}
We show that the \say{right} part of Theorem \ref{MainThm} holds. The left part is analogous. To prove Condition (1), we can, according to Condition (b), choose $k$ to equal $\pi(i)$ and hence satisfy Condition (i).

To prove (2), let us fix $i\in I$ and let $N^I_{\pi(i)}$ be the subset of $B_{i,{\pi(i)}}$ consisting of all elements $m$ such that $mB_{\pi(i), j}=0$ for all $ \pi(i)\neq j \in I$. Clearly, the simple right $R_{\pi(i)}$-module $T_{\pi(i)}$ is contained in $N^I_{\pi(i)}$.  Recall that Condition (a) on $R$ implies that $B_{i, \pi(i)}$ has a simple right socle; thus, the socle of $N^I_{\pi(i)}$ is also simple. 

The essentiality of the right socle of $eRe$ follows from Theorem \ref{ThmCombEss}, applying the fact that Conditions (a) and (b) hold in $R$. \end{dukaz}
 PF rings whose Nakayama permutation is a cycle were described in [14, Thm. 4]. In particular, it was shown that rings $R_1, \dots, R_n$ can be local corners of a ring whose Nakayama permutation is $(1~2~\dots~n)$ if and only if they form a \textit{cycle of dualities}: that is, if for each $i\leq n$, there exists a duality between left (linearly compact) $R_{i+1}$-modules and right (linearly compact) $R_i$-modules (modulo $n$). If one of the rings $R_i$ has a self-duality (e.g., it is a PF ring or an Artin algebra, such as a finite-dimensional algebra over a field or a finite ring), then all rings $R_i$ are isomorphic. It seems to be an open question whether a cycle of dualities can consist of non-isomorphic rings, even though local rings (even artinian) with a duality that does not possess self-duality exist [22].

A consequence of Müller's characterisation of the local rings on the diagonal is that all residue skew fields of local corners are isomorphic as fields. This is not explicitly stated in Müller's work. An alternative proof is in [9, Prop 11]. It is stated for QF rings, but it also applies in the PF case. We now generalise this observation further. Note that the endomorphism ring of the simple right $R$-module $T_i$ is isomorphic to the residue skew field of the local corner $R_i$.
\begin{Prop}\label{PropK}
    Let $R$ be a basic formal matrix ring of order $n$ with a Nakayama permutation $\pi$ and let $i\leq n$.

    If the right and left socles of $R$ coincide, then there exists a skew field $K$ such that the residue fields of local corners $R$ in the cycle of $\pi$ induced by $i$ are isomorphic to $K$.    
\end{Prop}

\begin{dukaz}
    Consider the simple socle $T:=T'_i=T_{\pi(i)}$ of $B_{i,\pi(i)}$ and denote by $K'$ the residue skew field of $R_i$ and by $K$ the residue skew field of $R_{\pi(i)}$. We prove that $K'\cong K$ as fields. The full statement of the proposition then follows by an iteration of this argument. 

    Because $T$ is simple on both sides, the projections of $R_i$ and $R_{\pi(i)}$ onto their residue skew fields induce a $K'\text-K$-bimodule structure on $T$ and $T$ is simple as a left $K'$-module and simple as a right $K$-module. 
    
    Fixing any nonzero element $m\in T$, we have $K'm=T=mK$. Because $T$ is faithful on both sides, this induces a bijection $\phi_m\colon K'\to K$ mapping $k\in K'$ to the unique element $\phi_m(k)\in K$ such that $km=m\phi_m(k)$. This bijection is a ring homomorphism and is hence an isomorphism of fields.  Indeed, take $a,b\in K'$ then \[abm=a(bm)=a(m\phi_m(b))=(am)\phi_m(b)=(m\phi_m(a))\phi_m(b)=m(\phi_m(a)\phi_m(b)).\]
\end{dukaz}

For the remainder of this subsection, let $R$ be a basic ring with a  Nakayama permutation $\pi=(1~2~\dots~n)$ whose right socle is an essential right ideal. To simplify the notation, when working with the indices $i, j, k$ of the formal matrix rings, we identify indices that are congruent modulo $n$

We have already discussed rings on the diagonal and pairing coordinate bimodules. As seen in Example \ref{ExCycle}, the remaining entries might be zero.  We show that if one bimodule is nonzero, then the whole corresponding shifted diagonal consists of nonzero bimodules, and often one more shifted diagonal does too. For PF rings, the relation of bimodules sharing the same pair of shifted diagonals can be further described in terms of duality, see [14, Thm. 4].
\begin{Prop}\label{PropCyclestructure}
    Let $R$ as above, $i, k\in \Z_n$ be indices such that $0\neq k\neq 1$ and assume that $B_{i, i+k}$ is nonzero. 

    Then for any $j\in \Z_n$, the bimodules $B_{j, j+k}$ and $B_{j, j-k+1}$ are nonzero.
\end{Prop}
If $2k\equiv 1$ in $\Z_n$, then the bimodules $B_{j, j+k}$ and $B_{j, j-k+1}$ coincide.
\begin{dukaz}
    Take a nonzero element $m\in B_{i, i+k}$. By Condition (b), $mB_{i+k, i+1}$ is nonzero. In particular,  $B_{i+k, i+1}$ must be nonzero. Applying (b) again to a nonzero element of $mB_{i+k, i+1}$, we find that $B_{i+1, i+k+1}$ is nonzero. The conclusion then follows from repeated use of (b). 
\end{dukaz}
As far as the support of the formal matrix ring is concerned, the above proposition cannot be strengthened any further. Any possibility that it allows can happen.
\begin{Prop}
    Let $n$ be a positive integer and let $I\subseteq \{2,3,\dots, n-1\}$ be a (possibly empty) index set.

    Then there exists a basic formal matrix ring $R$ with a Nakayama permutation $(1~2~\dots~n)$ and an essential right socle such that $B_{i,j}$ is nonzero if and only if it is either a local corner, a pairing coordinate bimodule, or $j\equiv i+k$ or $j\equiv i-k+1$ for some $k\in I$.
\end{Prop}
\begin{dukaz}
We construct such a ring as follows. Let $S$ be a local ring with self-duality defined by $S\text-S$-bimodule $E$ and let $\m$ be the residue skew field of $S$. Recall that $\m$ is an $(S\text-S)$-bimodule via the standard projection $S\twoheadrightarrow S/J(S)$ and both $_S\m$ and $\m_S$ are simple modules. Thus, we have an $(S\text-S)$-bimodule monomorphisms $\m\hookrightarrow E$.

We set $R_i:=S$ and $B_{i,i+1}:=E$ for any $i\in \Z_n$. For any $k\in I$, we set $B_{j, j+k}:=\m$ and $B_{j, j-k+1}:=\m$. We define all proper product homomorphisms as zero except for the pairing product homomorphisms $\varphi_{i, i+k, i+1}$ and $\varphi_{i,i-k+1, i+1}$. To define these, let $a,b\in \m$ and define $ab\hookrightarrow E$ where $ab$ is interpreted as multiplication in $\m$. 

To check the associativity condition ($\alpha$), observe that the multiplication of any three elements is always zero unless one of them is from the diagonal. The rest then follows from combinatorial criteria. 
\end{dukaz}
Although this proposition characterises all possible supports, the structure of product homomorphisms can differ. This example is based on  [20, Prop. 8.27, Thm. 6.15].
\begin{Ex}
    Let $Q$ be a cyclic quiver on the vertices $1, 2, 3$ with arrows $3\to2\to1\to3$ and let $K$ be a field. Consider the path algebra $KQ$ with the admissible ideal $I$ generated by paths of length at least six. 
    
Then $KQ/I$ is a basic Frobenius ring with a Nakayama permutation $(1~2~3)$. Representing it as a formal matrix ring where the $i$-th row corresponds to an indecomposable module generated by all paths starting at $i$, we see that $B_{i,i+1}B_{i+1,i}\neq 0$ and $B_{i,i+1}B_{i+1,i+2}\neq 0$.
\end{Ex}
It turns out that the assumption that the right socle is essential is indeed needed in Corollary \ref{CorFaithful} and Propositions \ref{PropTriangular} and \ref{PropCyclestructure}. 
\begin{Ex}
    Let $S$ be a local commutative ring that is not perfect and $E$ be the injective envelope of its simple module and $M$ be an $S$-module with zero socle.

    Then the trivial formal matrix rings \begin{gather*}
    \begin{pmatrix}
        S&E&M\\
        0&S&E\\
        E\oplus M&0&S
    \end{pmatrix},~~~~~~ \begin{pmatrix}
        S& E&0&0\\
        E&S&0&M\\
        0&0&S&E\\
        0&0&E&S
    \end{pmatrix}\end{gather*} are rings with Nakayama permutations $(1~2~3)$ and $(1~2)(3~4)$, respectively.
\end{Ex}

\section{Glueing of rings with respect to their Nakayama permutations}\label{SecGlue}

Under mild assumptions, this section demonstrates how to glue two indecomposable rings with a Nakayama permutation into a new indecomposable ring whose Nakayama permutation is the concatenation of their respective Nakayama permutations.

The method is as follows: first, we construct a block matrix ring with $S$ and $S'$ on the diagonal and copies of compatible skew fields in the off-diagonal blocks. This is described in Proposition \ref{PropGlue}. Because a ring with a Nakayama permutation decomposes as a block matrix ring along the cycles of the permutation, we can further strengthen this proposition: it is sufficient for the assumptions to be satisfied in only one cycle. This is formulated in Corollary \ref{CorGlue}. The main theorem of this section can be viewed as a special case of this corollary.

\smallskip

The first subsection provides a precise statement of Theorem \ref{ThmGlue} and discusses its assumptions and main applications. The construction is presented in the following subsections. Subsection \ref{SecFields} discusses compatible fields—a key assumption of the main theorem. The global construction is formulated as Proposition \ref{PropGlue} in Subsection \ref{SecCon}, and the generalisation is described in Subsection \ref{SecEndRings}. The main Theorem is used to construct an indecomposable Frobenius ring such that the endomorphism rings of simple modules are not isomorphic, see Example \ref{Exnon-isomfields}.

\subsection{Main Theorem and its aplications}\label{SecGlueThm}

We introduce the following auxiliary definitions, which will be used throughout this section. Examples of compatible skew fields are discussed in Subsection \ref{SecFields}.
\begin{Def}
        Two skew fields $K$ and $L$ are \emph{compatible} if $K$ admits a $K\text-L$-bimodule structure and $L$ admits an $L\text-K$-bimodule structure.

       Two semiperfect rings $S$ and $S'$ are \emph{locally compatible} if there exists a pair of compatible skew fields $K$ and $L$ such that $K$ is isomorphic to the top of a local corner of $S$ and $L$ is isomorphic to the top of a local corner of $S'$.
\end{Def}
  \begin{Def}
       Let $\sigma\in S_n$ and $\sigma'\in S_{n'}$ be two permutations. We can view them as permutations in $S_{n+n'}$ by letting them act as the identity on the extra elements. Let $\tau_n\in S_{n+n'}$ be a translation defined by $i\mapsto i+n \pmod{n+n'}$.

    Then the \emph{concatenation of $\sigma$ and $\sigma'$} is defined as $\pi:=\sigma \circ ( \tau_n \circ  \sigma' \circ \tau_n^{-1})\in S_{n+n'}$.
  \end{Def} 

The main result of this section is the following.
\begin{Thm}\label{ThmGlue}
Let $S$ and $S'$ be two indecomposable, non-simple, locally compatible semiperfect rings with Nakayama permutations $\sigma$ and $\sigma'$, respectively. Further, suppose that $S$ and $S'$ have essential socles.

Then there exists a semiperfect ring $R$ with an essential right socle whose Nakayama permutation is the concatenation of $\sigma$ and $\sigma'$, and there exists an idempotent $f\in R$ such that $fRf\cong S$ and $(1-f)R(1-f)\cong S'$.   
\end{Thm}
The assumption that $S$ and $S'$ are not simple is necessary for $R$ to be indecomposable, as shown in Proposition \ref{PropFixedPoint}. 

Ring $R$ has the same set of local corners as $S$ and $S'$. Recalling Proposition \ref{PropStructure}, $R$ is semiprimary or right perfect if and only if both $S$ and $S'$ are. This is not true for artinian rings, as all coordinate bimodules have to be finitely generated, as in Example \ref{ExFields2}.
\begin{Cor}
    The ring $R$ is semiprimary or right perfect if and only if both $S$ and $S'$ are.

    If $S$ and $S'$ are QF rings and $K_L$ and $L_K$ are finitely generated, then $R$ can be constructed as a QF ring.
\end{Cor}
By Example \ref{ExCycle}, given a field $K$ and a positive integer $n$, there exists a (Frobenius) ring whose Nakayama permutation is an $n$-cycle.
\begin{Cor}\label{CorPerm}
    Let $n$ be a positive integer and $\pi\in S_n$.
    
     Then, there exists an indecomposable (Frobenius) ring with a Nakayama permutation $\pi$.
\end{Cor}    
Building on the results of [8, Section 2], we proved that in an indecomposable right and left noetherian semiperfect ring, all endomorphism rings of simple modules have the same characteristic and are either all finite or have the same infinite cardinality [Já, Thm. 9]. If $R$ has an essential right socle, then the endomorphism rings of simple modules corresponding to the same cycle in the Nakayama permutation are isomorphic, by Proposition \ref{PropK}. 

We use Proposition \ref{PropGlue} to show that simple modules over a Frobenius ring need not have the same endomorphism ring if the Nakayama permutation does not link them.
\begin{Ex}\label{Exnon-isomfields}
    Let $L$ and $K$ be the non-isomorphic fields from Example \ref{ExFields2}. Let $S$ be a non-simple local Frobenius ring whose residue field is isomorphic to $K$ and $S'$ a non-simple local Frobenius ring whose residue field is isomorphic to $L$.

    Then $\begin{pmatrix}
        S&K\\
        L&S'
    \end{pmatrix}$ is a Frobenius ring whose Nakayama permutation is the identity. 
\end{Ex}
We conclude this subsection by observing that any skew field $K$ can be a residue field of a non-simple local Frobenius ring.  Consider its  \textit{trivial extension} by itself, denoted as $K\times K$, defined as the subring of the formal matrix ring $\begin{pmatrix}
    K & K\\

    0 & K
\end{pmatrix}$ consisting of matrices whose entries on the diagonal are equal. 

It is a local ring, and its residue field is isomorphic to $K$. Its Jacobson radical $E_{1,2} K$ is simple, which shows that it is a Frobenius ring.

\subsection{Examples of compatible fields}\label{SecFields}

Two isomorphic skew fields $K$ and $L$ are always compatible. It also suffices if $K$ contains an isomorphic copy of $L$ as a subring, and $L$ contains an isomorphic copy of $K$ as a subring. We present two such constructions.  We start with a folklore example from commutative algebra.

\begin{Ex}\label{ExFields1}
    Let $t_0, t_1, t_2,\dots$ be a countably infinite set of complex numbers algebraically independent over $\Q$. That is, they are transcendental over $\Q$ and do not satisfy any nontrivial polynomial equation with coefficients in $\Q$.

Let $K'$ be the algebraic closure of $\Q(t_0, t_1,\dots)$ in $\C$ and let $K$ be the algebraic closure of $\Q( t_1, t_2, \dots)$ in $\C$. Then we have field extensions $K\leq K(t_0) \leq K'$ with $K\cong K'$. However,  $K(t_0)$ is not isomorphic to $K$ because, by our assumption, $t_0$ is not algebraic over $K$; thus $K(t_0)$ is not algebraically closed.     
\end{Ex}

Note that in this case, $K(t_0)$ is not finitely generated and, therefore, is not linearly compact over $K$. We present a finite-dimensional construction due to R. Swan [21].

\begin{Ex}\label{ExFields2}
    Let $K'=\Q(x_1, \dots, x_{47})$ and $L$ be its subfield consisting of polynomials invariant under translation, i.e., $p(x_1, \dots, x_{47})\in K'$ is an element of $L$ if and only if  $p(x_1, \dots, x_{47})=p(x_{1+i}, \dots, x_{47+i \pmod{47})})$ for any $i<47$. 

    Let $K$ be a subfield of $L$ consisting of all symmetric polynomials, i.e., polynomials invariant under any $\pi\in S_{47}$. As an extension of $\Q$, it is generated by the elementary symmetric polynomials. Hence, both $K$ and $K'$ are transcendental extensions of $\Q$ of transcendental degree 47, and are therefore isomorphic. 

    By [21, Thm. 1], the field $L$ is not purely transcendental over $\Q$, and hence is not isomorphic to $K\cong K'$. 
    
    The field $K'$ is a finite-dimensional algebraic extension of $K$ as shown by the symmetric polynomial $\prod_{\substack{1\leq i \leq 47}}(x_i-y)$. From inclusions $K\leq L\leq K'$, it follows that $L$ is a finite-dimensional extension over $K$ and $K'$ is a finite-dimensional extension over $L$.
\end{Ex}

\subsection{Construction}\label{SecCon}

\begin{Prop}\label{PropGlue}
    Let $S$ and $S'$ be basic formal matrix rings of orders $n$ and $n'$ with Nakayama permutations $\sigma$ and $\sigma'$, respectively. Let $\pi$ be the concatenation of $\sigma$ and $\sigma'$.  Further, assume that 
    
    (A) All residue fields of local corner rings of $S$ are isomorphic to a fixed skew field $K$.

    (A') All residue fields of local corner rings of $S'$ are isomorphic to a fixed skew field $L$
    
    (B) The skew fields $K$ and $L$ are compatible.

   (C) For any $i\leq n$,  the simple submodules $T_{\sigma(i)}$ and $T'_i$ of pairing coordinate bimodules of $S$, defined in Conditions (2) and (2') of Theorem \ref{MainThm}, coincide. In addition, $T'_i$ is isomorphic to $K$ as an $S_i\text-S_{\sigma(i)}$-bimodule.

   (C') For any $i\leq n'$, the simple submodules $T_{\sigma'(i)}$ and $T'_i$ of pairing coordinate bimodules of $S'$, defined in Conditions (2) and (2') of Theorem \ref{MainThm}, coincide. In addition, $T'_i$ is isomorphic to $L$ as a $S'_i\text-S'_{\sigma'(i)}$-bimodule.

   (D) For any fixed point $i$ of $\sigma$, the local corner $S_i$ is not simple. 

      (D') For any fixed point $i$ of $\sigma'$, the local corner $S'_i$ is not simple. 
   
    Then there exists an indecomposable basic ring $R$ of order $n+n'$ with a Nakayama permutation $\pi$ and an idempotent $f\in R$ such that $fRf\cong S$ and $(1-f)R(1-f)\cong S'$.
\end{Prop}
Assumptions (D) and (D') are necessary to ensure that $R$ can be indecomposable; recall Proposition \ref{PropFixedPoint}.

\begin{Ex}Any two basic Frobenius $K$-algebras that are not division algebras can be glued together, cf. [14, Prop. 7].     
\end{Ex}

We construct ring $R$ as a block matrix ring of order two $\begin{pmatrix}
    S & K^{n\times n'}   \\
    L^{n'\times n} & S'
\end{pmatrix}$, where $K^{n\times n'}$ and $ L^{n'\times n}$ are blocks whose entries consist of copies of bimodules $_KK_L$ and $_LL_K$, respectively. 

Product homomorphisms within the blocks $S$ and $S'$ are thus inherited from the original rings. All \textit{new} proper product homomorphisms are set to zero, except for the pairing product homomorphisms of the forms $K^{n\times n'} L^{n'\times n}$ and $L^{n'\times n}K^{n\times n'}$.

Firstly, we need to establish how to view copies of $_KK_L$ and $_LL_K$ as bimodules over the corresponding local corners.
\begin{Lemma}
    Let $i\leq n < j\leq n+n'$.

    By Assumption (B), $K$ admits a $K\text-L$-bimodule structure. By Assumptions (A) and (A'), this induces an $S_i\text-S'_j$-bimodule structure on $K$ via the canonical projections $S_i\twoheadrightarrow \m_i\cong K$ and $S'_j\twoheadrightarrow \m'_j\cong L$.

 Similarly, $L$ admits an $S'_j\text-S_i$-bimodule structure.  
\end{Lemma}
We now move to defining pairing product homomorphisms  $\varphi_{i,j,\pi(i)}$, for $i\leq n <j\leq n+n'$. To this end, we use Assumption (C), first defining suitable bimodule homomorphisms to $K$. The remaining pairing product homomorphisms $\varphi_{j,i,\pi(j)}$, will be defined in an analogous fashion.
\begin{Lemma}\label{BimodLemma}
   Let $i, k\leq n<j\leq n+n' $.
   
   Then there exists an $S_i\text-S_k$-bimodule homomorphism $\Psi_{i,j,k}\colon  K \otimes_{S'_j} L \to K$ such that $a\otimes b$ is mapped to zero if and only if $a=0$ or $b=0$.
\end{Lemma}
We are interested only in the case $k=\pi(i)$, but the general proof works the same.
\begin{dukaz}
     Let $a\otimes_{S'_j} b\in   K \otimes_{S'_j} L$ be a pure tensor element. Then there exists $s\in S'_j$ such that the left $S'_j$ action of $s$ on $L$ maps $1_L$ to $b$. We can write
\begin{gather*}   
a\otimes_{S'_j} b=a\otimes_{S'_j} s\cdot 1_L=a\cdot s\otimes_{S'_j} 1_L
\end{gather*}
We define a nonzero abelian group homomorphism $K \otimes_{S'_j} L\to K$ by prescribing  $a\otimes_{S'_j}b\mapsto as$ and extending it additively to $K \otimes_{S'_j} L$. 

Note that $s$ from the preceeding equation is not defined uniquely, but its image under the standard projection $S'_j\twoheadrightarrow \m'_j$ is, and so is its right action on $K$.  This map is easily verified to be a homomorphism of $S_i\text-S_k$-bimodules.
\end{dukaz}

For $i\leq n < j\leq n+n'$, we define the $S_i\text-S_{\pi(i)}$-bimodule homomorphism \[\varphi_{i,j,\pi(i)}\colon _{S_i} K\otimes_{S'_j} L_{S_{\pi(i)}} \to B_{i,\pi(i)}\] as the composition of $\Psi_{i,j,\pi(i)}\colon K\otimes_{S'_j} L \to K$ and the injection $K\hookrightarrow B_{i,\pi(i)}$ from Assumption (C). 

Analogously, we define $\varphi_{j,i,\pi(j)}$ as  $S'_j\text-S'_{\pi(j)}$-bimodule homomorphism $L\otimes_{S_i}K\to L$ composed with $L\hookrightarrow B_{j,\pi(j)}$ from Assumption (C').

\smallskip

We now sketch the proof that $R$ is a well-defined formal matrix ring and satisfies the conclusion of Proposition \ref{PropGlue} using the combinatorial criterion.
\begin{dukaz}
\textbf{$R$ is a well-defined formal matrix ring. }

This amounts to showing the associativity condition $(\alpha)$. We address the cases \say{KLK} and \say{SKL} in detail. The rest is either inherited from rings $S$ and $S'$ or is analogous. 

(KLK) Let $i, k\leq n < j, l\leq n+n'$. Consider $a\in B_{i,j}\cong {_KK_L}$, $b\in B_{j,k}\cong {_LL_K}$, and $c\in B_{k,l}\cong {_KK_L}$. 

If $k\neq \pi(i)$, then $\varphi_{i,j,k}$ facilitates a product of the new off-diagonal blocks, which is not one of the defined pairing product homomorphisms and is thus set to zero, proving $(ab)c=0$. 

If $k=\pi(i)$, then $ab$ is in the pairing coordinate bimodule $B_{i,\pi(i)}$. The product $B_{i,k}B_{k, l}$ is defined to be zero as long as $k\neq i$. 

If $\pi(i)=i=k$, i.e., the local corner $S_i$ is a pairing coordinate bimodule, then $ab\in S_i$, and by Lemma \ref{LemmaJacob}, it is in the Jacobson radical of $S_i$. Now $B_{i,l}\cong {_KK_L}$ is simple as a left $S_i$-module, so any element of $J(S_i)$ annihilates it, showing $(ab)c=0$. 

The proof that $a(bc)=0$ is analogous.

 (SKL) Let $i, j, l\leq n <k\leq  n+n'$. Consider $a\in B_{i,j}$, $b\in B_{j,k}\cong {_KK_L}$, and $c\in B_{k,l}\cong {_LL_K}$.   If $i=j$, then the conclusion follows from Remark \ref{AssRemark}. Otherwise, $\varphi_{i,j,k}$ was defined as a zero morphism, so $(ab)c=0$.
 
 To show that $a(bc)=0$, first observe that if $l\neq \pi(j)$, then  $bc=0$. If $l=\pi(j)$, i.e., $\varphi_{j,k,l}$ is a newly defined pairing product homomorphism, then $bc$  is in the simple left $S_j$-submodule $T'_j$ of $B_{j,\pi(j)}$. By Condition (2') for $S$, we have $B_{i,j}(bc)=0$, so $a(bc)=0$.

 \smallskip 

\noindent \textbf{By Theorem \ref{MainThm}, $\pi$ is a Nakayama permutation of $R$.}

We show Conditions (1) and (2).  Conditions (1') and (2') are analogous.

    Let $j\neq \pi(i)$. If $B_{i,j}$ is a coordinate bimodule from the diagonal block, Condition (1) is inherited from the rings $S$ or $S'$. For the remaining coordinate bimodules, we can choose $k:=\pi(i)$. Condition (1) then follows directly  from Lemma \ref{BimodLemma}, since $\varphi_{i,j,\pi(i)}$ maps $a\otimes b$ to zero if and only if $a=0$ or $b=0$, by
    
   We prove Condition (2) for $i\leq n$. The case $i>n$ is analogous.
 By the construction of $\pi$, the coordinate bimodule $B_{i,\pi(i)}$ is in a diagonal block $S$.
 
 Condition (2) for $S$ ensures that the right $S_{\pi(i)}$-submodule of $B_{i,\pi(i)}$ consisting of elements that annihilate $B_{\pi(i),j}$ for any $\pi(i)\neq j\leq n$ has a simple socle $T_{\pi(i)}$. It remains to show that $T_{\pi(i)}B_{\pi(i),j}=0$ for $n<j\leq n+n'$.

If $i\neq \pi(i)$, then $\varphi_{i,\pi(i), j}$ is a zero homomorphism for any $j>n$, and (2) follows. If $i=\pi(i)$, then $B_{i,\pi(i)}=S_i$. The Jacobson radical of $S_i$ annihilates $_{S_i}K$, so it is enough to show that $T$ is contained in $J(S_i)$. Because $S_i$ is local, this holds whenever $J(S_i)$ is nonzero, which is equivalent to assuming $S_i$ is non-simple, as stated in (D).

\smallskip 

\noindent \textbf{$R$ is indecomposable.}

It is enough to test that the sums of canonical idempotents are not central. This follows because the off-diagonal blocks of $R$ contain only nonzero coordinate bimodules.
%\smallskip \noindent If $S$ and $S'$ have essential right socles, so does $R$.Following Theorem \ref{ThmCombEss}, condition (a) is inherited from the rings $S$ and $S'$. Condition (b) follows from our choice of pairing product homomorphisms and Lemma \ref{BimodLemma}
\end{dukaz}

    \subsection{Generalisation and proof of the main theorem}\label{SecEndRings}
    
We further strengthen Proposition \ref{PropGlue}, and subsequently prove the main theorem of this section.
\begin{Cor}\label{CorGlue}
    Let $S$ and $S'$ be two non-simple, indecomposable basic formal matrix rings of orders $n$ and $n'$ with Nakayama permutations $\sigma$ and $\sigma'$, respectively.  Let $\pi$ be the concatenation of $\sigma$ and $\sigma'$.

    Assume that there exist idempotents $e\in S$ and $e'\in S'$ such that  $eSe$ and $e'S'e'$ are rings with Nakayama permutations, which are the restrictions of $\sigma$ and $\sigma'$, respectively. Furthermore, assume that $eSe$ and $e'S'e'$ satisfy the assumptions of Proposition \ref{PropGlue}.

 Then there exists an indecomposable basic ring $R$ of order $n+n'$ with a Nakayama permutation $\pi$ and an idempotent $f\in R$ such that $fRf\cong S$ and $(1-f)R(1-f)\cong S'$. Furthermore, $R$ has an essential right socle if and only if $S$ and $S'$ have essential right socles. 
\end{Cor}
\begin{dukaz}
    Without loss of generality, we can assume that there exist $i\leq n$ and $j\leq n'$ such that $e=e_{n-i+1}+e_{n-i+2}+\dots+e_n$ and $e'=e'_1+e'_2+\dots +e'_{j}$.

Because $eSe$ and $e'S'e'$ satisfy assumptions of Proposition \ref{PropGlue}, we can construct a block matrix ring $\begin{pmatrix}
    eSe & K^{i\times j} \\
    L^{j\times i} & e'S'e'
\end{pmatrix}$ satisfying the conclusions of the said theorem. The ring $R$ is then constructed as a block matrix ring 

$\begin{pmatrix}
    (1-e)S(1-e) & (1-e)Se & 0 & 0\\
    eS(1-e) & eSe & K^{i\times j} &0 \\
    0 & L^{j\times i} & e'S'e' & e'S'(1-e')\\
    0&0& (1-e')S'e'& (1-e')S'(1-e')
\end{pmatrix}.$
The new products of form $(1-e)Se\otimes K^{i\times j}, K^{i\times j}\otimes e'S'(1-e'), L^{j\times i}\otimes eS(1-e)$, and $(1-e')S'e'\otimes L^{j\times i}$ are set to zero. 

The claim about the essentiality of socles follows easily from the combinatorial criterion, as the pairding product homomorphisms satisfy Condition (b) by definition. The rest is inherited from the rings $S$ and $S'$.
\end{dukaz}
We end this section with the proof of Theorem \ref{ThmGlue}.
\begin{dukaz} We can assume that $S$ and $S'$ are basic; otherwise, we can glue their basic rings into ring $R'$ and construct the formal matrix ring Morita equivalent to $R$ with multiplicities given by the multiplicities in the rings $S$ and $S'$.

We show that $S$ and $S'$ satisfy the assumptions of Corollary \ref{CorGlue}. Because they are locally compatible, there are indices $i, j$ such that $top(S_i)\cong K$ and $top(S'_j)\cong L$, where $K$ and $L$ are compatible skew fields. Let $e\in S$ and $e'\in S'$ be the sums of canonical idempotents of $S$ and $S'$ corresponding to the cycles of $\sigma$  and $\sigma'$ containing $i$ and $j$, respectively. 

Because $S$ and $S'$ are assumed to have an essential socle, by Theorem \ref{ThmCycle}, $eSe$ and $e'S'e'$ are rings whose Nakayama permutation is a cycle. By Proposition \ref{PropK}, all tops of local corners of $eSe$ and $e'S'e'$ are isomorphic to the skew fields $K$ and $L$, respectively. Thus $eSe$ and $e'S'e'$ satisfy Conditions (A) and (A') of Proposition \ref{PropGlue}. Condition (B) is satisfied by the assumptions on $K$ and $L$. Conditions (C) and (C') by Remark \ref{RemarkTT'}, and conditions (D) and (D') by Proposition \ref{PropFixedPoint}, since $S$ and $S'$ are assumed to be indecomposable non-simple.
\end{dukaz}

\medskip

\noindent\textbf{Aknowledgements:}
 The author wishes to express his thanks and appreciation to Jan Žemlička for his guidance and valuable insights during the research. 

\medskip

\noindent \textbf{Bibliography}

[1] Y. Baba, K. Oshiro. (2009). Classical Artinian Rings and Related Topics. (World Scientific Publishing). https://doi.org/10.1142/7451 

[2] Chen J, N. Ding, M. Yousif. (2004). On Noetherian Rings with Essential Socle. \textit{J. Aust. Math. Soc.} 76(1). 39-50. https://doi.org/10.1017/S1446788700008685

[3] C. Faith, D. Van Huynh. (2002).  When Self-injective Rings are QF: a Report on a Problem. 	\textit{J. Algebra Its Appl.} 1(1). 75-105. https://doi.org/10.1142/S0219498802000070

[4] K. Fuller. (1969). 
On indecomposable Injectives over Artinian Rings. Pac. J. Math. 29(1). 115-135 https://doi.org/10.2140/pjm.1969.29.115

[5] E. L. Green. (1982). On the Representation Theory of Rings in Matrix Form.  \textit{Pac. J. Math.} 100(1). 123-138. https://doi.org/10.2140/pjm.1982.100.123

[6] C. R. Hajarnavis, N. C. Norton. (1985). On Dual Rings and their Modules. \textit{J. Alg.} 93(2). 253-266. https://doi.org/10.1016/0021-8693(85)90159-0

[7] A. T. Hannula. (1973). On the Construction of Quasi-Frobenius Rings. \textit{J Alg.} 25(3). 403-414. https://doi.org/10.1016/0021-8693(73)90089-6

[8] M. Iovanov. (2016).   Frobenius–Artin algebras and Infinite Linear Codes. \textit{J. Pure Appl. Algebra}. 220(1). 560-576. https://doi.org/10.1016/j.jpaa.2015.05.030

[9] D. Krasula. (2026). Endomorphism rings of Simple Modules and Block Decomposition. \textit{J. Algebra Its Appl.}  Online Ready https://doi.org/10.1142/S0219498826501562

[10] P. Krylov, A. Tuganbaev. (2017). Formal Matrices. Algebra and Applications 23. (Springer). https://doi.org/10.1007/978-3-319-53907-2

[11] T. Lam. (1991). A First Course in Noncommutative Rings.  Graduate Texts in Mathematics 131. (Springer). https://doi.org/10.1007/978-1-4419-8616-0

[12] T. Lam. (1999). Lectures on Modules and Rings.  Graduate Texts in Mathematics 189. (Springer). https://doi.org/10.1007/978-1-4612-0525-8

[13] B. J. Müller. (1970). On Semi-perfect Rings.  \textit{Illinois J. Math} 14(3). 464-467. https://doi.org/10.1215/ijm/1256053082

[14] B. J. Müller. (1974). The Structure of Quasi-Frobenius Rings. \textit{	Can. J. Math.} 26(5). 1141-1151. https://doi.org/10.4153/CJM-1974-106-7

[15] T. Nakayama. (1941). On Frobeniusean Algebras II. \textit{Ann. Math.} 42(1). 1-21. https://doi.org/\\10.2307/1968984

[16] W. Nicholson, M. Yousif. (2003). Quasi-Frobenius Rings. 
Cambridge Tracts in Mathematics 158. (Cambridge University Press) 
https://doi.org/10.1017/CBO9780511546525

[17] I. Reiten. (1976). Stable Equivalence of Self-injective Algebras. \textit{J. Alg.} 40(1). 63-74. https://doi.org/\\10.1016/0021-8693(76)90087-9

[18] L. Rowen. (1988). Ring Theory I. Pure and Applied Mathematics 127. (Elsevier).

[19] B. Roux. (1971). Sur la Dualité de Morita. \textit{Tôhoku Math. J.} 23(3). 457-472. https://doi.org/10.2748/tmj/1178242594

[20] A. Skowroński, K. Yamagata. (2011). Frobenius Algebras I. EMS Textbooks in Mathematics. (EMS press). https://doi.org/10.4171/102

[21] R. Swan. (1969). Invariant Rational Functions and a Problem of Steenrod. \textit{Invent. Math.} 7. 148-158. https://doi.org/10.1007/BF01389798 

[22] W. Xue. (1989). Two Examples of Local Artinian Rings. \textit{Proc. Am. Math. Soc.} 107(1). 63-65. https://doi.org/10.2307/2048035

[23] W. Xue. (1992). Rings with Morita Duality. Lecture Notes in Mathematics 1523. (Springer) https://doi.org/10.1007/BFb0089071

[24] W. Xue. (1997). On a Theorem of Fuller. \textit{J. Pure Appl. Algebra}. 122(1-2). 159-168. https://doi.org/\\10.1016/S0022-4049(96)00070-9

\end{document}